\documentclass[12pt,twoside,leqno]{article}
\usepackage{amssymb}
\usepackage{enumerate}
\usepackage{amsmath,amsthm}
\newtheorem {thm}{Theorem}

\newtheorem*{prob}{Problem}

\theoremstyle{definition}
\newtheorem*{rem}{Remark}

\def\th{{^{\rm th}}}
 \def\D{{\partial}}
 \def\BCH{{\rm BCH}}
 \def\ad{{\rm ad}}
 \def\half{{\frac{1}{2}}}
 \def\ints{{\Bbb Z}}
 \def\lra{{\longrightarrow}}
 \def\takes{{\colon}}

\begin{document}
\textwidth=30cc \baselineskip=16pt

\title{A formula for topology/deformations and its significance}
\author{Ruth Lawrence\\
Einstein Institute of Mathematics, Hebrew University of Jerusalem\\
Givat Ram, Jerusalem 91904 ISRAEL\\
E-mail: ruthel@ma.huji.ac.il
 \and Dennis Sullivan\\
SUNY Department of Mathematics \& CUNY Graduate Center\\
Stony Brook, NY 11794-3651 USA \& New York NY 10016 USA\\
 E-mail: dennis@math.sunysb.edu}

\maketitle

\renewcommand{\thefootnote}{}


\renewcommand{\thefootnote}{\arabic{footnote}}
\setcounter{footnote}{0}
 %%%%%%%%%%%%%

\begin{abstract}
\noindent The formula is
$\D{e}=(\ad_e)b+\sum_{i=0}^\infty{\frac{B_i}{i!}}(\ad_e)^i(b-a)\>,$
with $\D{a}+{1\over2}[a,a] =0$ and $\D{b}+{1\over2}[b,b] =0$, where
$a$, $b$ and $e$ in degrees $-1$, $-1$ and 0 are the free generators
of a completed free graded Lie algebra $L[a,b,e]$. The coefficients
are defined by ${x\over{e^x-1}}=\sum_{n=0}^\infty{B_n\over{}n!}x^n$.
The theorem is that
\begin{itemize}
\item{this formula for $\D$ on
generators extends to a derivation of square zero on $L[a,b,e]$;}
\item{the formula for $\D{e}$ is unique satisfying the first property,
once given the formulae for $\D{a}$ and $\D{b}$, along with the
condition that the ``flow" generated by $e$ moves $a$ to $b$ in unit time.}
\end{itemize}
The immediate significance of this formula is that it computes the
infinity cocommutative coalgebra structure on the chains of the
closed interval. It may be derived  and proved using the geometrical
idea of flat connections and one parameter groups or flows of gauge
transformations. The deeper significance of such general DGLAs which
want to combine deformation theory and rational homotopy theory is
proposed as a research problem.
 \end{abstract}

\section{Introduction}

This paper fits into the general framework of constructions of
algebraic models of cell complexes using differential graded Lie
algebras and, conversely, topological and algebraic
interpretations of DGLAs. We give some background.

Let $L$ be a free Lie algebra on a graded vector space over {\bf Q}
provided with a derivation of degree plus one (or of degree minus
one).
%Note that an arbitrary differential Lie algebra $A$ can be
%viewed as a quotient of a free Lie algebra $L$, and the map
%$L\longrightarrow{}A$ induces an isomorphism from the homology of
%$L/[L,L]$ to the homology of $A/[A,A]$.

The choice of a differential of degree $-1$ with $L$ concentrated in
non-negative degrees is natural for the interpretation of these
objects in topology. If X is a cell complex with one 0-cell and only
2-, 3-, 4-, $\ldots$ cells, the rational homotopy theory of Quillen
\cite{Q} assigns a free differential Lie algebra $L$ with one
generator in degree $k$ for each $k+1$-cell ($k>0$). The homology of
$L$ is the Whitehead Lie algebra of homotopy groups tensor {\bf Q}
shifted down by one degree. The homology of $L/[L,L]$ is the
ordinary reduced homology ({\bf Q} coefficients) of the space
shifted down by one degree. One imagines that enlarging this
discussion to allow cells in degree 1 would be related to some Lie
algebras associated to non-trivial fundamental groups, but little is
known here, to our knowledge.

The choice of a differential of degree $+1$ with $L$ concentrated
in non-negative degrees is natural for the interpretation of such
differential Lie algebras controlling the deformation theory of
some mathematical structure. In this case one considers elements
in degree $1$ satisfying $dx+{1\over2}[x,x] =0$. One also makes
sense of the expression $x'= dy.y^{-1} + yxy^{-1}$ (gauge
transformation) for $y$ in (some completion of) degree $0$ and
declares $y$ to be an equivalence between $x$ and $x'$. The set of
equivalence classes is a formal version of the moduli space of the
structure whose deformations are controlled by $L$. One imagines
that enlarging the discussion with elements in degrees
$-1,-2,\ldots$ would involve degree $-1$ elements acting as
equivalences between equivalences, degree $-2$ elements as
equivalences between equivalences between equivalences etc, but
little is known to our knowledge.

The geometric or topological interpretation of a general
differential Lie algebra is a mixed object which combines both of
the above discussions: homotopy theory of spaces and moduli spaces
of deformations of some structure.

One knows that free Lie algebras arise from a standard (``bar")
construction starting from any differential cocommutative and
coassociative coalgebra over {\bf Q}. In this instance the
generators are those of the coalgebra shifted down by one. The
differential is the original differential extended to a derivation
plus the original comultiplication extended to a derivation. The
square of this derivation being zero is equivalent to
coassociativity. Then one may extend the notion of a cocommutative
coassociative coalgebra on a graded space to be any derivation of
square zero on the free Lie algebra (possibly completed) generated
by the graded space shifted down by one. This defines the notion
of an infinity differential graded cocommutative, coassociative
coalgebra. The higher terms beyond quadratic of the differential
are chain homotopies restoring coassociativity up to homotopy to the coproduct determined by the quadratic form.

Now one also knows that for any chain complex the cellular
approximations to the diagonal are homotopic, any two homotopies
between two of them are themselves homotopic, etc.  It follows
that there is on the chains an infinity cocommutative,
coassociative coalgebra structure.

\begin{prob}{Study this free differential Lie algebra
attached to a cell complex, determine its topological and geometric
meaning as an intrinsic object. Give closed form formulae for the
differential and for the induced maps associated to subdivisions.}
\end{prob}

In this note we will say something about the interval, the circle
and the real line. We use the Maurer-Cartan idea to find an explicit
formula involving Bernoulli numbers. We will see that the
subdivision map corresponding to splitting an interval into two by
adding a midpoint, is described by the Baker-Campbell-Hausdorff
formula.

\vskip1ex \noindent{\bf Previous work:} The abstract picture about
the diagonal goes back to Steenrod's construction of cohomology
operations mod $p$. The infinity coalgebra story was known for a
long time by Hinich et al. In the appendix to \cite{TZ} there is a
cell by cell canonical construction which is not explicit. In
\cite{CG}, there is an explicit sum over trees construction based on
Chen's iterated integrals and Whitney forms but not a closed form
expression.

\section{Preliminaries on flows on pre-DGLAs}

{\bf General DGLAs.} Recall that a DGLA is a vector space
$A$ over a field $k$ with grading $A=\oplus_{n\in\ints}A_n$ along
with a bilinear map $[.,.]\takes{}A\times{}A\lra{}A$ (bracket) and a
linear map $\D\takes{}A\lra{}A$ (differential) for which $\D^2=0$
while
 \begin{itemize}
 \item{(1)} (symmetry of bracket) $[b,a]=-(-1)^{|a||b|}[a,b]$;
 \item{(2)} (Jacobi identity)
$[[a,b],c]=[a,[b,c]]-(-1)^{|a||b|}[b,[a,c]]$;
 \item{(3)} (Leibniz rule) $\D[a,b]=[\D{a},b]+(-1)^{|a|}[a,\D{b}]$.
 \end{itemize}

Note that the three properties above are valid only for
homogeneous elements $a$, $b$ and $c$ of $A$ (that is an element of
$\bigcup_nA_n$), and $|a|\in\ints$ denotes the grading. The bracket
and differential are required to respect the grading, in that for
homogeneous elements, $|[a,b]|=|a|+|b|$ while $|\D{a}|=|a|-1$.  The
adjoint action of $A$ on itself is given by $\ad_e(a)=[e,a]$ and
acts on the grading by $\ad_e\takes{}A_n\lra{}A_{n+|e|}$. In this
notation (2) and (3) can be rewritten as
 \begin{itemize}
 \item{($2'$)} (Jacobi identity)
$\ad_{[a,b]}=[\ad_a,\ad_b]$
 \item{($3'$)} (Leibniz rule) $[\D,\ad_a]=ad_{\D{a}}$
 \end{itemize}

\noindent in which the brackets on the right-hand side refer to
the (signed) commutator of operators defined by
$[x,y]=xy-(-1)^{|x||y|}yx$ where the product is composition of
operators and the grading $|x|$ of a (homogeneous) operator is
the shift in grading which $x$ induces.  Thus the gradings of
$\ad_e$ and $\D$ are $|e|$ and $-1$, respectively.

When the condition $\D^2=0$ is removed, the resulting algebraic
structure will be called a {\it pre-DGLA}.

\vskip 1ex \noindent {\bf Auxiliary spaces over Q$[t]$.} For
simplicity we will work over $k={\bf Q}$, though the discussion also
holds for any field of characteristic 0. Assume now that $A$ is a
free Lie algebra on a finite dimensional graded vector space $V$, so
it has generators $x_1,\ldots,x_k$ where $x_i$ label a basis for
$V$. In order to deal with convergence issues which otherwise would
arise, we will need to work in certain finite dimensional quotients
of $A$. There is an additional grading on $A$ by the number of Lie
brackets,
$$A=\bigoplus_{n=0}^\infty A^{(n)}$$ in which $A^{(n)}$ is the
finite dimensional vector space generated by expressions involving
exactly $n$ Lie brackets in elements of $V$, so that it is spanned
by all words of the form
$$[x_{i_0},[x_{i_1},\ldots,[x_{i_{n-1}},x_{i_n}]\ldots]]$$ where
$i_0,\ldots,i_n\in\{1,\ldots,k\}$ label (not necessarily distinct)
basis elements of $V$, and there may be relations between them
induced by relations (1), (2). The grading is well-defined since
it is respected by relations (1), (2) and $A$ is free as a Lie
algebra.

For non-negative integers $N$, set
$B^{(N)}=\bigoplus_{n=0}^N{A^{(n)}}$; it has a natural Lie algebra
structure induced from $A$, and as such can be identified as the
(graded) Lie algebra quotient in which the vanishing of all
expressions involving exactly $N+1$ Lie brackets are imposed as
relations (as a consequence, all expressions with more than $N$
brackets must also vanish).  Then we have a tower of Lie algebra
homomorphisms
$$A\longrightarrow{}B^{(N)}\longrightarrow{}B^{(N-1)}\longrightarrow
\cdots\longrightarrow{}B^{(0)}=A^{(0)}=V\>,$$ with
$B^{(N)}\longrightarrow{}B^{(N-1)}$ mapping all elements of
$A^{(N)}$ to zero.

Define $U^{(N)}=\bigoplus_{n=0}^N(A^{(n)}\otimes{\bf Q}[t])$.
Picking a basis $\{{\bf e}_{n,r}|1\leq{}r\leq{}\dim{A^{(n)}}\}$ for
each $A^{(n)}$,  a typical element $x\in{}U^{(N)}$ can be written as
$$x=\sum_{n=0}^N\sum_{r=1}^{\dim{A^{(n)}}}p_{n,r}(t){\bf
e}_{n,r}$$ for some polynomials $p_{n,r}(t)\in{\bf Q}[t]$. Since
such an element involves only a finite number of such polynomials,
one can equivalently think of elements of $U^{(N)}$ as
 $$x=\sum_{m=0}^\infty{}t^m{\bf x}_m$$
 where only a finite number of the vectors ${\bf x}_m\in{}B^{(N)}$ are
 non-zero.  That is, an element of $U^{(N)}$ is a formal polynomial
 in $t$ with coefficients in $B^{(N)}$.

There is an obvious linear operator of differentiation by $t$
defined on $U^{(N)}$ by
$${d\over{}dt}\left(\sum_{m=0}^M{}t^m.{\bf
x}_m\right)=\sum_{m=0}^{M-1}t^m.(m+1){\bf x}_{m+1}\>.$$

\vskip 1ex \noindent{\bf  Differential structures over a free Lie
algebra.} For arbitrary elements $v_1,\ldots,v_k\in{}A$, there is
defined a unique linear map $\partial\takes{}A\longrightarrow{}A$
satisfying $\partial(x_i)=v_i$ for all $i$, along with the Leibniz
rule, that is, giving $A$ the structure of a pre-DGLA. The condition
that this defines the structure of a DGLA on $A$ is that
$\D\circ\D=0$, that is, that $\D^2(x)=0$ for all $x\in{}A$. Applying
the Leibniz rule twice gives $\D^2[u,v]=[\D^2u,v]+[u,\D^2v]$, from
which inductively it follows that a sufficient condition on $v_i$
for them to generate a DGLA structure on $A$ is that $\D^2x_i=0$ for
all $i$, that is $\D{v_i}=0$ for all $i$.

From freeness of $A$, it can inductively deduced that any such
pre-DGLA structure has
$\D(A^{(n)})\subset\bigoplus_{m\geq{}n}A^{(m)}$, so that it also
induces a well-defined pre-DGLA structure on $B^{(N)}$ which will be
a full DGLA structure so long as $\D^2$ vanishes on $A$ to order at
least $N$ in the Lie bracket grading.

\vskip 1ex \noindent{\bf Flatness and flows in a pre-DGLA.} For any
$v\in{}A_0$, consider the ``formal'' ordinary differential equation
 $$\frac{du}{dt}=\D{v}-\ad_v(u)\>,$$
for $u\in{}U^{(N)}$, where both sides are considered as elements
of $U^{(N)}$.    Writing $u(t)$ in the form
$\sum_{n=0}^\infty(t^n.x_n)$ where $x_n\in{}B^{(N)}$, the
differential equation breaks into the recurrence relation
\begin{eqnarray*}
 (n+1)x_{n+1}&=&-\ad_v(x_n),\qquad{}n>0\\
 x_1&=&\D{v}-\ad_v(x_0)
 \end{eqnarray*}
 from which we see that ($n\geq1$)
 $$x_n={(-\ad_v)^{n-1}\over{}n!}x_1={(-\ad_v)^{n-1}\over{}n!}(\D{v})
                        +{(-\ad_v)^n\over{}n!}x_0\>,$$
giving a unique solution for $u(t)\in{}U^{(N)}$ once the initial
condition $x_0=u(0)\in{}B^{(N)}$ is fixed.  Note that
$x_n=0\in{}B^{(N)}$ for all $n>N+1$, so that indeed $u(t)$ is {\it
polynomial} in $t$, for every choice of initial condition.

The differential equation, and hence also the
solution space, is invariant under time-translation and so is said
to define the {\sl flow} on $U^{(N)}$ generated by $v$. By
evaluation at a given $t=t_0\in{\bf Q}$, this flow defines an
action of $({\bf Q},+)$ on $U^{(N)}$ by
  $$t_0.x_0\equiv{}u(t_0)\>,$$
namely the action of $t_0$ on an element $x_0\in{}U^{(N)}$ is
given by ``flowing according to the flow generated by $v$ for time
$t_0$''.  Observe that despite the intuition based on a continuous
model and derivatives, the formal definitions are only valid at
{\it rational} ``times'' $t$, and that this works because of the
rationality of coefficients in all the expansions. Explicitly,
$$t_0.x_0=x_0+\sum_{n=1}^{N+1}t_0^n\left({(-\ad_v)^{n-1}\over{}n!}(\D{v})
                        +{(-\ad_v)^n\over{}n!}x_0\right)
$$

An element $x\in{}A_{-1}$ is said to be {\sl flat} iff
$\D{x}+\half[x,x]=0$. Similarly, if this equality holds up to
order $N$ brackets, then it is a flat element of $B^{(N)}_{-1}$.

\section{The DGLA model for the interval}

The interval consists of two points and a single 1-cell.  Its model
should therefore have two generators in degree $-1$ (corresponding
to its endpoints) and a single generator in degree 0.

\begin{thm} There is a unique completed free
differential graded Lie algebra, $A$, with generating elements $a$,
$b$ and $e$, in degrees $-1$, $-1$ and $0$ respectively, for which
$a$ and $b$ are flat while the flow generated by $e$ moves from $a$
to $b$ in unit time. The differential is specified by,
$$\D{e}=(\ad_e)b+\sum_{i=0}^\infty{\frac{B_i}{i!}}(\ad_e)^i(b-a)\>,$$
where $B_i$ denotes the $i\th$ Bernoulli number defined as
coefficients in the expansion
${x\over{e^x-1}}=\sum\limits_{n=0}^\infty{}B_n {x^n\over{}n!}$.
\end{thm}

\begin{proof}

Let $A$ be the free graded Lie algebra generated by $a$, $b$ and
$e$. For any non-negative integer $N$, define the derived spaces
$B^{(N)}$ and $U^{(N)}$ as in the previous section. We prove the
result of the theorem on the truncated free Lie algebra $B^{(N)}$
for all $N$ and see that the corresponding differentials are
compatible for all $N$.  For any $x\in{}B^{(N)}_{-1}$, define a
map $\D_x\takes{}A\lra{}A$ by its action on the generators
 $$\D_x{a}=-\half[a,a],\qquad\D_x{b}=-\half[b,b],\qquad\D_x{e}=x\>,$$
extended to the whole of $A$ via linearity and the Leibniz rule.
(Here $\D_x$ is well-defined since derivation by the Leibniz rule
preserves relations (1) and (2).) This defines a pre-DGLA structure
on $B^{(N)}$.

The flow on $B^{(N)}$ generated by $e\in{U^{(N)}_{-1}}$ has
 $\frac{du}{dt}=x-\ad_e(u)$. For the particular solution with $u(0)=a$, the solution is given
as in the previous section by
$$u(t)=a+\sum_{n=1}^{N+1}t^n\left({(-\ad_e)^{n-1}\over{}n!}x
                        +{(-\ad_e)^n\over{}n!}a\right)
=e^{-t.\ad_e}a+{e^{-t.\ad_e}-1\over(-\ad_e)}x
$$
where as always the operator exponential is defined by its series
expansion $e^X=\sum\limits_{n=0}^\infty{X^n\over{}n!}$ so that the
operator acting on $x$ is defined by its series expansion
$\sum\limits_{n=1}^\infty{\frac{t^n}{n!}}(-\ad_e)^{n+1}$. Observe
that although there are formally infinite series in the expression
for $u(t)$, as elements of $U^{(N)}$ they are finite sums. The
condition that $u(1)=b$ is precisely that
$$x=\frac{(-\ad_e)}{e^{-\ad_e}-1}\big(b-e^{-\ad_e}a\big)
   =(\ad_e)b+\frac{ad_e}{e^{\ad_e}-1}(b-a)\>,$$
namely the value of $\D{e}$ given in the proposition (compatible
elements of $B^{(N)}$ for different $N$). This proves uniqueness.

It remains only to verify existence, that is to show that the
pre-DGLA structure is in fact a full DGLA structure, that is,
$\D_x^2=0$ for this particular value of $x$. From the Leibniz rule
for $\D_x$, it follows that for all $p,q\in{}B^{(N)}$,
$\D_x^2[p,q]=[\D_x^2p,q]+[p,\D_x^2q]$, so that it is only necessary
to check that $\D_x^2=0$ on the generators. For the generator $a$,
we have
$$\D_x^2(a)=\D_x(-\half[a,a])=[a,\D_xa]=[a,-\half[a,a]]=0\>,$$
the final equality following from the Jacobi identity. Similarly
$\D_x^2(b)=0$.

To prove that $\D_x^2(e)=0$, consider the flow $u$ generated by
$e$ for which $u(0)=a$ as above.  By our choice of $x$, this flow
also has $u(1)=b$. Consider the function $f(t)=\D_xu+\half[u,u]$
(the curvature), taking values in $B^{(N)}_{-2}$ at rational $t$;
equivalently, $f$ defines an element of $U^{(N)}_{-2}$. Its
derivative is
 \begin{eqnarray*}
  \frac{df}{dt}
  &=&\D_x\frac{du}{dt}+\big[u,\frac{du}{dt}\big]\\
  &=&\D_x(x-\ad_eu)+[u,x-\ad_eu]\\
  &=&\D_x^2(e)-\big(\ad_{\D_xe}(u)+\ad_e\D_xu\big)+([u,x]-[u,\ad_eu])\\
  &=&\D_x^2(e)-\ad_e(f(t))
 \end{eqnarray*}
    where we have used the Leibniz rule and that
$[u,\ad_eu]=\half\ad_e[u,u]$ from the Jacobi identity.  Thus $f$
satisfies a first-order linear differential equation with constant
(operator) coefficients of the same form as that satisfied by $u$
where now $x$ is replaced by $\D_x^2e$, while $f(0)=f(1)=0$ (since
$a$ and $b$ are flat). It follows that $\D_x^2e=0$ in $U^{(N)}_{-2}$
(for all $N$), as required.
 \end{proof}

\begin{rem} From the last calculation in the proof, it can be
seen that in the DGLA of Theorem 1, a flow $u$ defined by an
arbitrary element $v\in{}B^{(N)}_0$ on $B^{(N)}_{-1}$ has curvature
$f\in{}U^{(N)}_{-2}$ satisfying ${df\over{}dt}=-\ad_vf$, so that if
$u(0)$ is flat (so that $f(0)=0$) then $f$ is identically zero. That
is, the flow on $A_{-1}$ defined by an arbitrary element of $A_0$,
preserves flatness.
\end{rem}

\begin{rem} The vanishing of odd Bernoulli numbers after the first
is exactly the condition to make the formula for $\D{}e$ in terms of
$a$, $b$ display the symmetry of the interval, that is be anti-invariant under interchange of $a$, $b$ accompanied by a
sign change of $e$,
$$\D{e}=(b-a)+{1\over2}\ad_e(a+b)+{1\over12}(\ad_e)^2(b-a)-{1\over720}(\ad_e)^4(b-a)+\cdots\>.$$
 In other words, if a flow in direction
$e$ moves $a$ to $b$ in unit time, then $-e$ moves $b$ to $a$ in
unit time.
\end{rem}

\vskip 1ex \noindent{\bf Differential geometric interpretation.} We
would like to push the language of `flatness' and `flows' used in
the above formal proof a little further.  The real differential
geometric meaning of these constructions is however yet to be
understood.

The proof of Theorem 1 used the language of curvature and flatness
of connections, alluding to  the interpretation of the Maurer-Cartan
equation
$$\D{a}=-{1\over2}[a,a]$$ as the condition for a connection $a$
(as given by its associated 1-form) to be flat. However, there is
also a deformation theory interpretation, namely that the
differential structure in a DGLA $(A,{\D})$ can be deformed by
replacing $\D$ by $\D_a\equiv\D+\ad_a$ on condition that the
Maurer-Cartan condition is satisfied. Indeed, the deformed
structure defines a DGLA so long as $\ad_{\D_ax}=[\D_a,\ad_x]$ and
$\D_a^2=0$; we calculate
\begin{eqnarray*}
\ad_{\D{}x+[a,x]}
&=&ad_{\D{}x}+ad_{[a,x]}=[\D,\ad_x]+[\ad_a,\ad_x]=[\D+\ad_a,\ad_x]\\
(\D+\ad_a)^2
&=&\D^2+\D.\ad_a+\ad_a.\D+(\ad_a)^2=[\D,\ad_a]+ad_{{1\over2}[a,a]}=ad_{\D{a}+{1\over2}[a,a]}
\end{eqnarray*}
(by repeated applications of the Jacobi identity and the Leibniz
rule) so that the first condition is automatic, while the vanishing
of $\D_a^2$ is guaranteed by the Maurer-Cartan equation.

So we will think of any $a\in{}A_{-1}$ as defining a ``connection''
$\D_a=\D+\ad_a$. Furthermore, any $e\in{}A_0$ generates a flow on
$u\in{}A_{-1}$, which we consider as an ``infinitesimal gauge
transformation" flowing connections by ${du\over{}dt}=\D{}e-\ad_eu$;
as we saw from the proof of Theorem 1, this flow preserves flatness
of connections.  In addition, $e\in{}A_0$ also defines a flow on
$v\in{}A_0$ by ${dv\over{}dt}=-\ad_ev$, for which
\begin{eqnarray*}
{d\over{}dt}(\D_uv)
&=&{d\over{}dt}(\D{}v+\ad_uv)=\D{dv\over{}dt}+\ad_u{dv\over{}dt}+\ad_{du\over{}dt}v\\
&=&\D{}(-\ad_ev)+\ad_u(-\ad_ev)+\ad_{\D{}e}v-\ad_{\ad_e{}u}v\\
&=&-\ad_e\D{}v-\ad_e\ad_uv=-\ad_e(\D_u{}v)
\end{eqnarray*}
so that the condition $\D_uv=0$ is preserved by the flow, that is,
the flow $e$ defines a ``parallel transport" between fibres
$\ker\D_u\subset{}A_0$ over each ``point'' $u(t)$.

 \section{The DGLA model of the circle}

 To obtain a circle from an interval we need  only identify the
 endpoints.  So we obtain a single 0-cell $a$ and a single 1-cell
 $e$.  In our algebraic model, we impose the condition $a=b$ and
 immediately the differential of Theorem 1 collapses. The
 resulting model of the circle is a free Lie algebra with two
 generators $a$ and $e$ in degrees $-1$ and $0$, respectively,
 and differential $\D$:
 $$\D{a}=-{1\over2}[a,a]\>,\qquad\D{e}=\ad_ea\>.$$
The twisted differential over the point $a$ would then be just
$\D_a=\D+\ad_a$ and here $\D_ae=0$ so that the localization to a
point (namely to $a$) would be just generated by $e$ in grading 0
-- this corresponds to the single generating loop $e$.

So our model now corresponds points to flat connections; in
particular, the endpoints of the interval are corresponded to $a$
and $b$, while interior rational points $t_0$ give other flat
elements $u(t_0)$. On the other hand, the 1-cell is represented by
$e$ which defines an infinitesimal gauge transformation flowing
between all these (rational) points (flat connections).
Furthermore, as will be discussed in a further paper, this model
can be extended to higher dimensions, and then it will become
apparent that the algebraic analogue of localization to a point
$a$ (corresponding to consideration of the loop space based at
that point) is the replacement of the whole complex with
differential $\D$, by the complex truncated to non-negative
degrees with degree 0 part restricted to $\ker\D_a$ and
differential $\D_a$.
 A path connecting points then induces a flow preserving the complexes.

\section{Gluing intervals}

Suppose $X$ is any 1-complex.  Using Theorem 1, we can construct its
DGLA model as a free Lie algebra on generators $a_i$ in degree $-1$
for each 0-cell, and $e_i$ in degree 0 for each 1-cell, along with a
differential $\D$ which is uniquely defined by its action on the
generators making $a_i$ flat and giving $\D{}e_i$ by the
corresponding formula from Theorem~1 in which $a$ and $b$ are
replaced by the algebra elements associated with the endpoints of
the interval labelled by $e_i$.

For example, a 1-complex consisting of two adjoining intervals,
$[a_0,a_1]$ and $[a_1,a_2]$ with corresponding 1-cells $e_1$ and
$e_2$, is modelled by the DGLA $B$, which, as a Lie algebra, is
free on generators $a_0$, $a_1$, $a_2$ in degree $-1$ and $e_1$,
$e_2$ in degree 0. The differential $\D$ has
\begin{eqnarray*}
\D{}a_0
&=&-{1\over2}[a_0,a_0]\>,\qquad\D{}a_1=-{1\over2}[a_1,a_1]\>,\qquad\D{}a_2=-{1\over2}[a_2,a_2]\>,\\
\D{e_1}&=&(\ad_{e_1})a_1+\sum_{i=0}^\infty{\frac{B_i}{i!}}(\ad_{e_1})^i(a_1-a_0)\>,\\
\D{e_2}&=&(\ad_{e_2})a_2+\sum_{i=0}^\infty{\frac{B_i}{i!}}(\ad_{e_2})^i(a_2-a_1)\>.
\end{eqnarray*}

The geometric model is a subdivided interval. The previous
discussion defines a flow on $B_{-1}$; the flow according to $e_1$
flows from $a_0$ to $a_1$ in unit time and then the flow according
to $e_2$ flows from $a_1$ to $a_2$ in unit time.  That is, to any
rational point $p$ on either interval, there corresponds a flat
element $u_p\in{}A_{-1}$.

Removing the interior point $a_1$ would give a single interval
$[a_0,a_2]$ with 1-cell $e$ whose DGLA model, $A$, is a free Lie
algebra on generators $a_0$, $a_2$ in degree $-1$ and $e$ in
degree $0$, along with differential $\D$ defined by
\begin{eqnarray*}
\D{}a_0
&=&-{1\over2}[a_0,a_0]\>,\qquad\D{}a_2=-{1\over2}[a_2,a_2]\>,\\
\D{e}&=&(\ad_{e})a_2+\sum_{i=0}^\infty{\frac{B_i}{i!}}(\ad_{e})^i(a_2-a_0)\>.
\end{eqnarray*}
To parallel the geometric fact that the glued pair of intervals is
just a subdivision of a single interval, we have the following
theorem.

\begin{thm}There is homomorphism
$p\takes{}A\longrightarrow{}B$ respecting the DGLA structure, for
which $p(a_0)=a_0$, $p(a_2)=a_2$ while $p(e)$ is given by the
Baker-Campbell-Hausdorff formula on $e_1$ and $e_2$,
$$p(e)=e_1+e_2+{1\over2}[e_1,e_2]
+{1\over12}[e_1,[e_1,e_2]]-{1\over12}[e_2,[e_1,e_2]]+\cdots\>.$$
\end{thm}

\begin{proof} Denote by $\BCH(e_1,e_2)$ the
Baker-Campbell-Hausdorff formula on $e_1$ and $e_2$. By the Jacobi
identity, it follows that as operators,
$ad_{\BCH(e_1,e_2)}=\BCH(\ad_{e_1},\ad_{e_2})$ and hence that (again
as operators), $e^{\ad_{\BCH(e_1,e_2)}}=e^{\ad_{e_1}}e^{\ad_{e_2}}$.

By the Leibniz
rule, to prove that $p$ is a DGLA homomorphism, it is enough to
check compatibility of the action of the differential on
generators, that is, that $p(\D{}e)=\D(p(e))$. As we saw in the
proof of Theorem~1, $\D{}e_1$, $\D{}e_2$ satisfy
\begin{eqnarray*}
a_1
&=&e^{-\ad_{e_1}}a_0+{{e^{-\ad_{e_1}}-1}\over-\ad_{e_1}}(\D{e_1})\>,\\
a_2
&=&e^{-\ad_{e_2}}a_1+{{e^{-\ad_{e_2}}-1}\over-\ad_{e_2}}(\D{e_2})\>.
\end{eqnarray*}
Substituting the first equation into the second gives the identity
in $B$:
$$a_2=e^{-\ad_{e_2}}e^{-\ad_{e_1}}a_0
+e^{-\ad_{e_2}}{{e^{-\ad_{e_1}}-1}\over-\ad_{e_1}}(\D{e_1})
+{{e^{-\ad_{e_2}}-1}\over-\ad_{e_2}}(\D{e_2})\>.$$ On comparson with
the identity $a_2
=e^{-\ad_{e}}a_0+{{e^{-\ad_{e}}-1}\over-\ad_{e}}(\D{e})$ in $A$, and
recalling from above that
$e^{-ad_{p(e)}}=e^{-ad_{e_2}}e^{-ad_{e_1}}$, the theorem now follows
from Lemma~3 in the Appendix.
\end{proof}

\begin{rem} According to the previous discussion, $e_1$ induces a
flow on $B_{-1}$ (points/connections) and also on $B_0$. If we think
of $B_0$ as being in the fibre over the corresponding point, the
flow on $v\in{}B_0$ is defined by the differential equation
${dv\over{}dt}=-\ad_{e_1}v$. This is a homogeneous linear
differential equation, and so its solution is
$v(t)=e^{-t.\ad_{e_1}}v(0)$, so that the twisting in the fibre as we
move along the edge $e_1$ from $a_0$ to $a_1$ is $e^{-\ad_{e_1}}$
and a `flat section' has values $v_0$, $v_1$ over the endpoints of
the interval $[a_0,a_1]$ related by $v_1=e^{-\ad_{e_1}}v_0$.
Similarly $e_2$ induces a flow from $a_1$ to $a_2$ as well as on the
corresponding ``fibres" and in composition
$$v_2=e^{-\ad_{e_2}}v_1=e^{-\ad_{e_2}}e^{-\ad_{e_1}}v_0=e^{-ad_{p(e)}}v_0\>.$$
Thus Theorem~2 implies that not only is the map
$A\longrightarrow{}B$ a DGLA homomorphism, but it is also
compatible with the induced ``flow" structures discussed
previously.
\end{rem}

\section{Appendix: Some algebraic lemmas}

We here give the proofs for three lemmas used in the previous
section.

\vskip1ex\noindent{\bf Lemma 1.}\quad\emph{Let $e$ and $v$ be
arbitrary elements of a DGLA with $e$ of degree 0.  Then the
following formula holds, where $E$ denotes the operator $\ad_e$:
$$\D(e^{-E}v)=e^{-E}(\D{}v)+(-1)^{|v|}e^{-E}ad_v{e^{E}-1\over{}E}(\D{}e)\>.$$}

\begin{proof} By the definition of $e^{-E}$,
$$\D{(e^{-\ad_e}v)}
=\sum_{n=0}^\infty{(-1)^n\over{}n!}\D{\big({(\ad_e)}^nv\big)}\>.$$
However
\begin{align*}
&\D{\big({(\ad_e)}^nv\big)}\\
 &=\D([e,[e,\cdots[e,v]\cdots]])\\
 &=[\D{e},[e,\cdots[e,v]\cdots]]+[e,[\D{e},\cdots[e,v]\cdots]]+\cdots
    +[e,[e,\cdots[\D{e},v]\cdots]]\\
    &\qquad\qquad+[e,[e,\cdots[e,\D{v}]\cdots]]\\
 &=(-1)^{|v|+1}\bigg(\ad_{(\ad_e)^{n-1}v}(\D{e})
     +E.\ad_{(\ad_e)^{n-2}v}(\D{e})
     +\cdots+E^{n-1}\ad_v(\D{e})\bigg)\\
     &\qquad\qquad+E^n(\D{v})\>.
\end{align*}
By the Jacobi identity, the operator
$\ad_{(\ad_e)^mv}=\ad_{[e,\cdots[e,v]\cdots]}$ can be equivalently
written as the repeated commutator of operators
$$[\ad_e,\cdots[\ad_e,\ad_v]\cdots]=[E,\cdots[E,V]\cdots]$$ where
$V\equiv\ad_v$. Since $e$ is of degree zero, the grading on the
operator $E$ is also zero and so these graded operator commutators
are all with standard sign $[A,B]=AB-BA$. Inductively one obtains
$$\ad_{E^mv}=(\ad_E)^mV=\sum_{k=0}^m(-1)^k{m\choose{k}}E^{m-k}VE^k\>.$$
Substituting into the above expression
\begin{eqnarray*}
\D{\big({(\ad_e)}^nv\big)}
 &=&(-1)^{|v|+1}\bigg(\sum_{m=0}^{n-1}E^{n-m-1}\ad_{E^mv}\bigg)(\D{e})+E^n(\D{v})\\
 &=&(-1)^{|v|+1}\bigg(\sum_{m=0}^{n-1}\sum_{k=0}^m(-1)^k{m\choose{}k}E^{n-k-1}VE^k\bigg)(\D{e})+E^n(\D{v})\\
 &=&(-1)^{|v|+1}\bigg(\sum_{k=0}^{n-1}(-1)^k{n\choose{}k+1}E^{n-k-1}VE^k\bigg)(\D{e})+E^n(\D{v})\>.
\end{eqnarray*}
Combining over all $n$, this gives
\begin{eqnarray*}
\D{(e^{-\ad_e}v)}
 &=&\sum_{n=0}^\infty{(-1)^n\over{}n!}
 \bigg((-1)^{|v|+1}\bigg(\sum_{k=0}^{n-1}(-1)^k{n\choose{}k+1}E^{n-k-1}VE^k\bigg)(\D{e})
 +E^n(\D{v})\bigg)\\
 &=&(-1)^{|v|}\bigg(\sum_{n=1}^\infty\sum_{k=0}^{n-1}{(-1)^{n-k-1}\over{}n!}
    {n\choose{}k+1}E^{n-k-1}VE^k\bigg)(\D{e})+e^{-E}(\D{v})\\
 &=&(-1)^{|v|}e^{-E}V{e^E-1\over{}E}(\D{e})+e^{-E}(\D{v})\>.
\end{eqnarray*}\end{proof}

\vskip1ex\noindent{\bf Lemma 2.}\quad\emph{Let $e$ and $v$ be
arbitrary elements of a DGLA with $e$ of degree 0.  Then
$\ad_{(e^{-E}v)}=e^{-E}Ve^{E}$ where $E\equiv\ad_e$,
$V\equiv\ad_v$.}

\begin{proof}
 From the expression for $\ad_{E^mv}$
 obtained in the proof of Lemma~1,
\begin{eqnarray*}
\ad_{(e^{-E}v)}
 &=&\sum_{n=0}^\infty{(-1)^n\over{}n!}\ad_{E^nv}\\
 &=&\sum_{n=0}^\infty\sum_{k=0}^n{(-1)^{n+k}\over{}n!}{n\choose{k}}E^{n-k}VE^k\\
 &=&e^{-E}Ve^E\>.
\end{eqnarray*}\end{proof}

\vskip1ex\noindent{\bf Lemma 3.}\quad\emph{If $e_1$ and $e_2$ are
elements of an arbitrary DGLA with $f=\BCH(e_1,e_2)$, then the
following identity holds:
$${{e^{-\ad_{f}}-1}\over-\ad_{f}}(\D{f})
=e^{-\ad_{e_2}}{{e^{-\ad_{e_1}}-1}\over-\ad_{e_1}}(\D{e_1})
+{{e^{-\ad_{e_2}}-1}\over-\ad_{e_2}}(\D{e_2})\>.
$$}

\begin{proof} We will derive this identity as a compatibility condition. It
follows from Lemma~1, directly from the definition of the
exponential and the Leibniz rule, that for arbitrary elements $e$,
$v$ of a DGLA in which $e$ has degree $0$,
$$\D(e^{-E}v)=e^{-E}(\D{}v)+(-1)^{|v|}e^{-E}ad_v{e^{E}-1\over{}E}(\D{}e)\>,$$
where $E$ denotes the operator $\ad_e$. Applying this to evaluate
$\D(e^{-\ad_f}v)$ in two different ways, one directly from this
formula, and the other as $\D(e^{-\ad_{e_2}}e^{-\ad_{e_1}}v)$ by
applying the formula twice, and equating the results, gives that for
all Lie algebra elements $v$, there is an identity
$$e^{-F}\ad_v{e^F-1\over{}F}(\D{}f)
=e^{-E_2}e^{-E_1}\ad_v{e^{E_1}-1\over{}E_1}(\D{}e_1)
+e^{-E_2}\ad_{(e^{-E_1}v)}{e^{E_2}-1\over{}E_2}(\D{}e_2)\>,$$ where
again we have denoted the operators $\ad_{e_1}$, $\ad_{e_2}$ and
$\ad_f$ by $E_1$, $E_2$ and $F$, respectively. Furthermore, by
Lemma~2, $\ad_{(e^{-E_1}v)}=e^{-E_1}\ad_ve^{E_1}$, while
$e^{-F}=e^{-E_2}e^{-E_1}$ is an invertible operator, so that the
relation becomes
$$\ad_v{e^F-1\over{}F}(\D{}f)
=\ad_v{e^{E_1}-1\over{}E_1}(\D{}e_1)
+\ad_{v}e^{E_1}{e^{E_2}-1\over{}E_2}(\D{}e_2)\>.$$ Since this holds
in all DGLAs, it holds in particular for the free DGLA on three
generators $e_1$, $e_2$, $v$ (all in degree 0), and we can therefore
remove the $\ad_v$ operator from all terms, leaving an equality,
which when multiplied by $e^{-F}=e^{-E_2}e^{-E_1}$ on the left, is
the required identity.
\end{proof}

\subsection*{Acknowledgements}

The possibility of such a differential on the free Lie algebra on
three generators corresponding to the unit interval arose for one of
us almost twenty years ago in conversations with Maxim Kontsevich at
the IHES. At first it was considered to be a remarkable
combinatorial miracle involving Bernoulli numbers.  Later it was
regarded by several workers to be an obvious consequence of the
homotopy commutativity and associativity of the diagonal map for the
chains on the unit interval.

The authors worked out in Israel (2001) the proof that the formula
for $\D{e}$ that takes the flat element $a$ to the flat element $b$
by the flow associated to the infinitesimal gauge group action leads
to a derivation of square zero, but its publication was somewhat
delayed by the addition to our respective families of (three, one)
wonderful children! The first author would like to thank Stony Brook
for its hospitality during her visit in 2010, where the continuation
of the story from an interval through 1-complexes to more general
complexes was pursued. The authors would like to thank the referees
for useful comments.

\end{document}